\documentclass[a4paper,12pt]{article}

\usepackage{amssymb}
\usepackage{amsthm}
\usepackage[mathscr]{eucal}
\usepackage[english]{babel} 
\usepackage[tbtags]{amsmath}
\usepackage{enumerate}
\usepackage[matrix,arrow,tips,curve]{xy}

\usepackage[dvips]{graphics}
 \usepackage{psfrag}
 \psfrag{uno}{\fontsize{30}{30}$\Gamma_1$}
 \psfrag{due}{\fontsize{30}{30}$\Gamma_2$}
 \psfrag{tre}{\fontsize{40}{40}$\Gamma_Z$}
 \psfrag{C1}{\fontsize{20}{20}$C_1$}
 \psfrag{C2}{\fontsize{20}{20}$C_2$}
 \psfrag{C3}{\fontsize{20}{20}$C_3$}
 \psfrag{Z}{\fontsize{20}{20}$Z$}

\newtheorem{Theorem}{Theorem}
\newtheorem{Lemma}[Theorem]{Lemma}
\newtheorem{Proposition}[Theorem]{Proposition}
\newtheorem{Corollary}[Theorem]{Corollary}

\theoremstyle{definition}
\newtheorem{Definition}[Theorem]{Definition}

\theoremstyle{remark}

\newtheorem*{Example}{Example}

\newcommand{\Pic}{\operatorname{Pic}}

\newcommand{\Hom}{\operatorname{\emph{Hom}}}
\newcommand{\Sing}{\operatorname{Sing}}

\newcommand{\Aut}{\mathrm{Aut}}

\newcommand{\Hilb}{\mathrm{Hilb}}

\newcommand{\Z}{\mathbb{Z}}

\newcommand{\C}{\mathbb{C}}
\renewcommand{\P}{\mathbb{P}}

\newcommand{\mgo}{\mathcal{M}_g^0}

\newcommand{\mgbar}{\overline{\mathcal{M}}_g}

\newcommand{\pry}{\operatorname{\mathit{Pr}}}
\newcommand{\im}{\operatorname{Im}}
\newcommand{\pdg}{P_{d,g}}
\newcommand{\pdgbar}{\overline{P}_{d,g}}

\newcommand{\pdgfun}{\overline{\mathcal{P}}_{d,g}}
\newcommand{\prymbar}{\overline{\pry}_g}

\newcommand{\rg}{\overline{R}_g}

\newcommand{\prymx}{\pry_Z}
\newcommand{\prymfun}{\overline{\operatorname{\mathcal{P}r}}_g}
\newcommand{\ph}{\ensuremath{\varphi}}
\newcommand{\w}{\widetilde}
\newcommand{\pr}[1]{\ensuremath{\mathbb{P}^{#1}}}
\newcommand{\ol}[1]{\ensuremath{\mathcal{O}_{#1}}}

\title{Moduli of Prym curves}
\author{Edoardo Ballico, Cinzia Casagrande and 
Claudio Fontanari\thanks{This research is part of the 
T.A.S.C.A. project of I.N.d.A.M., supported by P.A.T. (Trento) and 
M.I.U.R. (Italy).}}
\date{}

\begin{document}
\maketitle

\begin{small}
\begin{center}
\textbf{Abstract}
\end{center}
Here we focus on $\rg$, the compactification of the moduli space of curves 
of genus $g$ together with an unramified double cover, constructed by 
Arnaud Beauville in order to compactify the Prym mapping.
We present an alternative description of $\rg$, inspired by the moduli 
space of spin curves of Maurizio Cornalba, and we discuss in detail 
its main features, both from a geometrical and a combinatorial point of view. 
\vspace{0.5cm}

\noindent
\textsc{A.M.S. Math. Subject Classification (2000)}: 14H10, 14H30.
\end{small}

\subsection*{Introduction}
Let $X$ be a smooth curve of genus $g$.
As it is well known (see for instance \cite{beau}, p.~104, or 
\cite{acgh}, Appendix B, \S~2, 13.), a square root of $\mathcal{O}_X$ 
corresponds to an unramified double cover of $X$.  

A compactification $\rg$ 
of the moduli space of curves of genus $g$ 
together with an unramified double cover was constructed by 
Arnaud Beauville (\cite{beauville}, Section~6; see also
\cite{prym},  Theorem~1.1) by means of admissible double covers of
stable curves. This moduli space was introduced
as a tool to compactify the mapping which associates to a curve 
plus a $2$-sheeted cover the corresponding Prym variety; 
however, we believe that it is interesting also in its own and 
worthy of a closer inspection.

Here we explore some of the geometrical and combinatorial 
properties of $\rg$. In order to do that, we present a 
description of this scheme 
which is different from the original one and is 
inspired by the construction performed by Maurizio Cornalba in 
\cite{cornalba1} of the moduli space of spin curves $\overline{S}_g$.
This is a natural compactification over $\mgbar$ of the space of
pairs $(X,\zeta)$, where $\zeta\in\Pic X$ is a square root of 
the canonical bundle $K_X$.

In Section~\ref{def} we define a \emph{Prym curve} 
to be just the analogue of a spin curve; Cornalba's approach 
can be easily adapted to our context and allows us to put a 
structure of projective variety on the set $\prymbar$ of 
isomorphism classes of Prym curves of genus $g$. 
This variety has two irreducible components $\prymbar^{\,-}$ and
$\prymbar^+$, where $\prymbar^{\,-}\simeq\mgbar$ contains 
``trivial'' Prym curves; moreover, by comparing Prym curves and admissible 
double covers, we give an explixit isomorphism between $\prymbar^+$ and
$\rg$ over $\mgbar$.
 
Next, in Section~\ref{geo} we reproduce the arguments in 
\cite{fontanari} in order to show that $\prymbar$ 
is endowed with a natural
injective morphism
into the compactification of the universal Picard 
variety constructed by Lucia Caporaso in \cite{caporaso}, just like
$\overline{S}_g$. 

Finally, in Section~\ref{com} we turn to the combinatorics of 
$\prymbar$. Applying the same approach used in \cite{spin} for 
spin curves, we study the ramification of the morphism
$\prymbar\rightarrow\mgbar$ over the boundary.
We describe
the numerical properties of the scheme-theoretical fiber $\prymx$ over
a point $[Z]\in\mgbar$, which turn out to 
depend only on the dual graph $\Gamma_Z$ of $Z$. 
From this combinatorial description, it follows that  
the morphisms $\prymbar\rightarrow\mgbar$ and
$\overline{S}_g\rightarrow\mgbar$ ramify in a different way.

The  moduli space $\rg$ of admissible double covers has
been studied also by Mira Bernstein in \cite{mira}, where $\rg$ is
shown to be of general type for $g=17,19,21,23$ (\cite{mira},
Corollary 3.1.7) (for
$g\geq 24$ it is obvious, since $\mgbar$ is).

We work over the field $\C$ of complex numbers.

We wish to thank Lucia Caporaso for many fruitful conversations.
We are also grateful to the anonymous referee for pointing out a gap
in a previous version of this paper.
%%%%%%%%%%%%%%%%%%%%%%%%%%%%%%%%%%%%%%%%%%%%%%%%%%%%%%%%%%%%%%%%%%%%%
\section{Prym curves and admissible double covers}\label{def}
{\bf \ref{def}.1. Defining the objects.}
Let $X$ be a Deligne-Mumford semistable curve and $E$ 
an irreducible component of $X$. One says that $E$ is 
\emph{exceptional} if it is smooth, rational, and meets the other 
components in exactly two points. Moreover, one says that $X$ is 
\emph{quasistable} if any two distinct exceptional components of $X$ 
are disjoint. The stable model of $X$ is the stable curve $Z$ obtained
from $X$ by contracting each exceptional component to a point.
In the sequel, $\w{X}$ will denote the subcurve 
$\overline{X \smallsetminus \cup_i E_i}$ obtained from $X$ by removing all  
exceptional components. 

We fix an integer $g\geq 2$.
\begin{Definition}
A \emph{Prym curve} of genus $g$
is the datum of $(X,\eta,\beta)$ where $X$ is 
 a quasistable curve of genus $g$, $\eta\in\Pic X$, and 
$\beta\colon \eta^{\otimes 2} \to \mathcal{O}_X$ is a 
homomorphism of invertible 
sheaves  satisfying the following conditions:
\begin{enumerate}[(i)]
\item $\eta$ has total degree $0$ on $X$ and degree $1$ on every exceptional 
component of $X$;
\item $\beta$ is non zero at a general point of every non-exceptional 
component of $X$.
\end{enumerate}
We say that $X$ is the support of the Prym curve $(X,\eta,\beta)$.

An \emph{isomorphism}
between two Prym curves $(X, \eta,\beta)$ and 
$(X', \eta',\beta')$ 
is an isomorphism
$\sigma\colon X \to X'$ 
such that there exists an isomorphism $\tau\colon \sigma^*(\eta') \to \eta$
 which makes the following diagram commute 
%$\sigma^*(\mathcal{O}_{X'})$ and $\mathcal{O}_X$.
\footnote{Observe that
  we are adopting the  
  convention  that the datum of $\tau$ is not
  included in the definition of isomorphism, as in \cite{cornalba1}. 
This is different from
  the convention in \cite{cornalba2}; see \cite{cornalba2}, end of
  section 2, for
  a discussion about this.} 
\[
\xymatrix{
{\sigma^*(\eta')^{\otimes 2}} \ar[d]_{\sigma^*(\beta')}
\ar[r]^(.65){\tau^{\otimes 2}}
&{\eta^{\otimes 2}}\ar[d]^{\beta}\\
{\sigma^*(\ol{X'})} \ar[r]^(.6){\sim}&{\ol{X}}.}
\]
\end{Definition}
Let $(X, \eta,\beta)$ be a Prym curve and let $E_1,\dotsc,E_r$ be
the exceptional components of $X$.
From the definition it follows that $\beta$ vanishes identically 
on all $E_i$ and induces an isomorphism 
$$
\eta^{\otimes 2} \vert_{\w{X}} 
\stackrel{\sim}{\longrightarrow}
\mathcal{O}_{\w{X}}(-q^1_1-q^2_1- \cdots -q^1_r-q^2_r),
$$ 
where $\{q^1_i,q^2_i\}=\w{X}\cap E_i$ for $i=1,\dotsc,r$.
In particular, when $X$ is smooth, $\eta$ is just a point of order two
in the Picard group of $X$. The number of such points, as it is well-known, 
is exactly $2^{2g}$.

We denote by $\Aut(X,\eta,\beta)$ the group of automorphisms of the Prym
curve $(X,\eta,\beta)$. 
As in \cite{cornalba1}, p. 565, one can show that
$\Aut(X,\eta,\beta)$ is finite. 

We say that an isomorphism between two Prym
curves $(X,\eta,\beta)$ and $(X,\eta',\beta')$
having the same support  is 
\emph{inessential} if it induces the identity on the stable model of
$X$. We denote by $\Aut_0(X,\eta,\beta)\subseteq\Aut(X,\eta,\beta)$ 
the subgroup
of inessential automorphisms. 
We have the following
\begin{Lemma}[\cite{cornalba1}, Lemma 2.1]
\label{inessential}
There exists an inessential isomorphism between two Prym curves
 $(X,\eta,\beta)$ and $(X,\eta',\beta')$  if and only if 
$$\eta\vert_{\w{X}}\simeq \eta'
\vert_{\w{X}}.$$ 
\end{Lemma}
\noindent So the set of isomorphism classes of
Prym curves supported on $X$ is in bijection with the set of square
roots of $\mathcal{O}_{\w{X}}(-q^1_1-q^2_1- \cdots -q^1_r-q^2_r)$ in
$\Pic\w{X}$, modulo the action of the group of automorphisms of 
$\w{X}$ fixing $q^1_1,q^2_1,\dotsc,q^1_r,q^2_r$.

A \emph{family of Prym curves} is a flat family of quasistable 
curves $f\colon\mathcal{X} \to S$ with an invertible sheaf
$\boldsymbol{\eta}$ on $\mathcal{X}$ and a homomorphism
$$
\boldsymbol{\beta}\colon \boldsymbol{\eta}^{\otimes 2} 
\longrightarrow \mathcal{O}_{\mathcal{X}}
$$
such that the restriction of these data to any fiber of $f$ gives rise 
to a Prym curve.
An isomorphism between two
  families of Prym curves $(\mathcal{X} \to
  S,\boldsymbol{\eta},\boldsymbol{\beta})$  and 
$(\mathcal{X}' \to S,\boldsymbol{\eta}',\boldsymbol{\beta}')$ over $S$
 is an isomorphism
$\sigma\colon \mathcal{X} \to \mathcal{X}'$ over $S$ such that there exists
an isomorphism
$\tau\colon \sigma^*(\boldsymbol{\eta'}) \to \boldsymbol{\eta}$
 compatible with the
canonical isomorphism between $\sigma^*(\mathcal{O}_{\mathcal{X}'})$
and $\mathcal{O}_{\mathcal{X}}$.

We define the moduli functor associated to Prym curves in the obvious
way:  $\prymfun$ is the contravariant functor from the category of
schemes to the one of sets, 
which to every scheme $S$ associates the set $\prymfun(S)$ of isomorphism 
classes of families of Prym curves of genus $g$ over $S$. 

\bigskip 

\noindent {\bf \ref{def}.2. The universal deformation.} 
Fix a Prym curve $(X, \eta,\beta)$, 
call $Z$ the stable model of $X$ and denote by 
$E_1, \dotsc, E_r$ the exceptional components of $X$. 
Let $\mathcal{Z}'\rightarrow B'$ 
be the universal deformation of $Z$, where
 $B'$ is the unit policylinder in $\C^{3g-3}$ with coordinates
$t_1,\dotsc,t_{3g-3}$ such that  $\{t_i=0\}\subset B'$ 
is the
locus where the node corresponding to $E_i$ persists for $i=1,\dotsc,r$.
Let $B$ be another unit policylinder in $\C^{3g-3}$ with coordinates
$\tau_1,\dotsc,\tau_{3g-3}$, and consider the map 
$B\rightarrow B'$ given  by
$t_i=\tau_i^2$ for $1\leq i\leq r$ and $t_i=\tau_i$ for
$i>r$. Call $\mathcal{Z}$ the pull-back of
$\mathcal{Z}'$ to $B$. For $i\in\{1,\dotsc,r\}$ the family 
$\mathcal{Z}_{|\{\tau_i=0\}}\rightarrow \{\tau_i=0\}\subset B$ has a
section $V_i$, 
 corresponding to the locus of the $i$th node. Let
 $\mathcal{X}\rightarrow \mathcal{Z}$ be the blow-up of
 $V_1,\dotsc,V_r$ and call $\mathcal{E}_1,\dotsc,\mathcal{E}_r$ the
 exceptional divisors.
\begin{equation*}
\xymatrix{
{\mathcal{X}}\ar[r]\ar[rd]&{\mathcal{Z}}\ar[d]\ar[r]&{\mathcal{Z}'}
\ar[d]
\\&{B}\ar[r]&{B'}}
\end{equation*}
The variety $\mathcal{X}$ is smooth and
$\mathcal{X}\rightarrow B$ is a family
of quasistable curves,  with $X$ as central fiber.
 Up to an inessential automorphism, we can assume that 
 $\eta^{\otimes
   2}\simeq\mathcal{O}_{\mathcal{X}}(-\sum_i\mathcal{E}_i)_{|X}$ and
that this isomorphism is induced by $\beta$. By 
  shrinking $B$ if necessary, we can extend $\eta$ to
$\boldsymbol{\eta}\in\Pic\mathcal{X}$ such that 
 $\boldsymbol{\eta}^{\otimes 2}
\simeq\mathcal{O}_{\mathcal{X}}(-\sum_i\mathcal{E}_i)$.
Denote by $\boldsymbol{\beta}$ the composition of this isomorphism
with the natural inclusion
$\mathcal{O}_{\mathcal{X}}(-\sum_i\mathcal{E}_i)\hookrightarrow
\mathcal{O}_{\mathcal{X}}$. Then $(\mathcal{X} \to B,
\boldsymbol{\eta},\boldsymbol{\beta})$  
is a family of Prym curves, and there is a
morphism $\psi\colon X\to\mathcal{X}$ which induces an isomorphism
of Prym curves between $(X,\eta,\beta)$ and the fiber of the family
over $b_0=(0,\dotsc,0)\in B$. 
This family 
provides a \emph{universal deformation} of $(X,\eta,\beta)$: 
\begin{Theorem}
\label{univdef}
Let $(\mathcal{X}'\to T,\boldsymbol{\eta}',\boldsymbol{\beta}')$ be a
family of Prym curves and let $\ph\colon X\to\mathcal{X}'$ be a morphism
which induces an isomorphism of Prym curves between $(X,\eta,\beta)$
and the fiber of the family over $t_0\in T$.

Then, possibly after shrinking $T$, there exists a unique morphism
$\gamma\colon T\to B$ satisfying the following conditions:
\begin{enumerate}[(i)]
\item $\gamma(t_0)=b_0$;
\item there is a cartesian diagram
$\qquad\xymatrix{
{\mathcal{X}'}\ar[d]\ar[r]^{\delta} & {\mathcal{X}}\ar[d] \\
T\ar[r]^{\gamma} & {B\,;}}$
\item $\boldsymbol{\eta}'\simeq \delta^*(\boldsymbol{\eta})$ and 
$\boldsymbol{\beta}'=\delta^*(\boldsymbol{\beta})$;
\item $\delta\circ\ph=\psi$.
\end{enumerate}
\end{Theorem}
\noindent Since the proof of \cite{cornalba1},
Proposition 4.6 applies verbatim to our case, we omit the proof of
Theorem \ref{univdef}.

\bigskip

\noindent {\bf \ref{def}.3. The moduli scheme.}
Let $\prymbar$ be the set of isomorphism classes of Prym curves of 
genus $g$. 
%and $\pry_g$ be the subset consisting of classes of smooth curves. 
We define a natural structure of analytic variety on $\prymbar$ 
following \cite{cornalba1}, \S~5.

Consider a Prym curve $(X,\eta,\beta)$ and its universal deformation
$(\mathcal{X}\to B,\boldsymbol{\eta},\boldsymbol{\beta})$ constructed
in 1.2. By the
universality, 
the group $\Aut(X,\eta,\beta)$ acts on $B$ and on $\mathcal{X}$. 
This action has the
following crucial property:
\begin{Lemma}[\cite{cornalba1}, Lemma 5.1]
\label{quoz}
Let $b_1,b_2\in B$ 
and let $(X_{b_1},\eta_{b_1},\beta_{b_1})$ and $(X_{b_2},
\eta_{b_2},\beta_{b_2})$  be the fibers of the universal
family over $b_1$ and $b_2$ respectively.
 Then
there exists $\sigma\in\Aut(X,\eta,\beta)$ such
that $\sigma(b_1)=b_2\,$\footnote{Where we still denote by $\sigma$ the
  automorphism of $B$ induced by $\sigma$.} if and only if
the Prym curves $(X_{b_1},\eta_{b_1},\beta_{b_1})$ and $(X_{b_2},
\eta_{b_2},\beta_{b_2})$ are isomorphic.
\end{Lemma}
\noindent Lemma \ref{quoz} implies that the natural 
(set-theoretical)
map $B\to\prymbar$, associating to $b\in B$ the isomorphism class of the
fiber over $b$, descends to a well-defined, injective map
$$J\colon B / \Aut(X, \eta,\beta)\ \longrightarrow\ \prymbar.$$
 This allows to define a complex
structure on the subset $\im J\subseteq\prymbar$. 
Since $\prymbar$ is covered by these subsets, in order to get a
complex structure on $\prymbar$ we just have to
check that the  complex structures are compatible on the overlaps.

This compatibility will follow from the following remark,
which is an immediate consequence of 
the construction of the universal family in 1.2:
\begin{enumerate}[$\bullet$]
\item
{\em the family of Prym curves
$(\mathcal{X}\to B,\boldsymbol{\eta},\boldsymbol{\beta})$ is a
universal deformation for any of its fibers.}
\end{enumerate}
In fact, assume that there are two Prym curves $(X_1,\eta_1,\beta_1)$
and $(X_2,\eta_2,\beta_2)$ such that the images of the associated maps
$J_1$, $J_2$ intersect. Choose a Prym curve $(X_3,\eta_3,\beta_3)$
corresponding to a point in $\im J_1\cap\im J_2$. Let $B_i$
($i=1,2,3$)
be the
basis of the universal deformation of $(X_i,\eta_i,\beta_i)$. Then by
the remark above, for $i=1,2$ 
 there are
natural open immersions $h_i\colon B_3\hookrightarrow B_i$, 
equivariant with respect to the
actions of the automorphism groups. Hence $h_i$ induces an open
immersion $\overline{h}_i\colon
B_3/\Aut(X_3,\eta_3,\beta_3)\hookrightarrow
B_i/\Aut(X_i,\eta_i,\beta_i)$, and $J_3=J_i\circ\overline{h}_i$.

Observe now that
the morphisms 
$$ B/\Aut(X,\eta,\beta)\longrightarrow B'/\Aut(Z)$$
glue together and yield a morphism $p\colon\prymbar\to\mgbar$.
Clearly $p$ is finite,
 as a morphism between analytic varieties (see \cite{SGA1EXPXII}). 
Hence $\prymbar$ is projective, because
$\mgbar$ is.
The variety $\prymbar$ has finite quotient singularities; in
particular, it is normal.

The degree of $p$ is
$2^{2g}$. The fiber over a smooth curve $Z$ is just the
set of points of order two in its Picard group, modulo the action of
$\Aut(Z)$ if non trivial. 
When $Z$ is a stable
curve, the set-theoretical fiber over $[Z]$ consists of isomorphism
classes of Prym curves $(X,\eta,\beta)$ such that the stable model of $X$ is
$Z$.  
In section \ref{com}
we will describe precisely the scheme-theoretical fiber over $[Z]$, 
following \cite{cornalba1} and \cite{spin}.
We will show that $p$ 
is \'etale over $\mgbar^0\smallsetminus D_{irr}$, where 
$\mgbar^0$ is the locus of stable curves with trivial automorphism
group and $D_{irr}$ is the boundary component whose general
member is an irreducible curve with one node.

Finally, $\prymbar$ is a coarse moduli space for the functor
$\prymfun$. For any family of Prym curves over a scheme $T$,
the associated moduli morphism $T\to \prymbar$ is locally defined by Theorem
\ref{univdef}.

Let $\prymbar^{\,-}$ be the closed subvariety of $\prymbar$ consisting
of classes of Prym curves $(X,\eta,\beta)$ where $\eta\simeq\ol{X}$.
 Observe that when
 $\eta$ is trivial, the curve $X$ is stable. So
$\prymbar^{\,-}$ is the image of the obvious section of
$p\colon\prymbar\to\mgbar$, and it is an irreducible (and connected)
component of $\prymbar$, isomorphic to $\mgbar$.

Let $\prymbar^{+}$ be the complement of $\prymbar^{\,-}$ in $\prymbar$, and
denote by $\pry_g^+$ its open subset consisting of
classes of Prym curves supported on smooth curves.
Then $\pry_g^+$ parametrizes connected
unramified double covers of 
smooth curves of genus $g$; it is well-known that this moduli space is
irreducible, being a finite quotient of the moduli space of
smooth curves of genus $g$ with a level 2 structure, which is
irreducible by \cite{dm}. So $\prymbar^{+}$ is an irreducible
component of $\prymbar$.

\bigskip

\noindent {\bf \ref{def}.4. Admissible double covers.}
Consider a pair 
$(C,i)$ where $C$ is a stable curve of genus $2g-1$ and $i$ is
an involution of $C$ such that:
\begin{enumerate}[$\bullet$]
\item the set $I$ of fixed points of $i$ is contained in $\Sing C$;
\item for any fixed node, $i$ does not exchange the two
  branches of the curve. 
\end{enumerate}
Then the quotient $Z:=C/i\,$ is a stable curve of genus $g$, and
$\pi\colon 
C\rightarrow Z$ is a finite morphism of degree 2,
 \'etale over $Z\smallsetminus \pi(I)$.
This is called an \emph{admissible double cover}.
Remark that $\pi$ is not a cover in the usual sense, since it is not flat 
at $I$.

The moduli space $\rg$ of admissible double covers of stable curves of
genus $g$ is constructed \cite{beauville}, Section 6 
(see also \cite{prym,prymhulek}), as the moduli space for pairs
$(C,i)$ as above.
 
An isomorphism of two admissible covers $\pi_1\colon 
C_1\rightarrow Z_1$ and $\pi_2\colon C_2\to Z_2$
 is an isomorphism $\ph\colon Z_1\stackrel{\sim}{\to}Z_2$ such
that there exists\footnote{Given $\ph$, there are exactly two choices for
  $\w{\ph}$; if $\w{\ph}$ is one, the other is
 $\w{\ph}\circ i_1$.} an isomorphism 
$\w{\ph}\colon C_1\stackrel{\sim}{\to} C_2$ with $\pi_2\circ\w{\ph}
=\ph\circ \pi_1$. 

We denote by $\Aut(C\rightarrow Z)$ the automorphism group of the
 admissible cover $C\to Z$, so
 $\Aut(C\rightarrow Z)\subseteq\Aut (Z)$. All elements of 
 $\Aut(C\rightarrow Z)$ are induced by automorphisms of
$C$, different
from $i$, and that commute with $i$.

Let $(C,i)$ be as above; we describe its universal deformation.
 Let $\mathcal{C}'\rightarrow W'$ be a universal
deformation of $C$. By the universality, there are compatible
involutions $i'$ of $W'$ and $\boldsymbol{i}'$ of
$\mathcal{C}'$, extending the action of $i$ on the central fiber.
 Let $W\subset
W'$ be the locus fixed by $i'$,  $\mathcal{C}\rightarrow W$ the induced
family and $\boldsymbol{i}$ the restriction of $\boldsymbol{i}'$ to
 $\mathcal{C}$. Then  
 $(\mathcal{C},\boldsymbol{i})\rightarrow W$ is a universal deformation of
$(C,i)$ and the corresponding family of admissible double covers is 
 $\mathcal{C}\rightarrow\mathcal{Q}:=\mathcal{C}/\boldsymbol{i}
\rightarrow W$.
 
We are going to show that $\rg$ is isomorphic over $\mgbar$
to the irreducible component $\prymbar^+$ of $\prymbar$.

First of all we define a map
$\Phi$ from the set of
non trivial Prym curves of genus $g$ to the set of 
admissible double covers of  stable curves of genus $g$.

Let $\xi=(X,\eta,\beta)$ be a Prym curve
with $\eta\not\simeq\mathcal{O}_X$; 
then $\Phi(\xi)$ will be an admissible double cover of the stable
model $Z$ of $X$, constructed as follows.
The homomorphism $\beta$ induces an isomorphism
$$\eta_{|\w{X}}^{\otimes (-2)}  
\simeq\mathcal{O}_{\w{X}}(q^1_1+q^2_1+ \cdots +q^1_r+q^2_r).$$
This determines a double cover
$\w{\pi}\colon
 \w{C}\rightarrow\w{X}$, ramified over $q^1_1,q^2_1,\dotsc,q^1_r,q^2_r$, which
are smooth points of $\w{X}$.   
Now call
$C_{\xi}$ 
the stable curve obtained identifying $\w{\pi}^{-1}(q^1_i)$ 
with $\w{\pi}^{-1}(q^2_i)$ for all
$i=1,\dotsc,r$. 
Then the induced map $C_{\xi}
\rightarrow Z$ is the admissible double cover $\Phi(\xi)$.

Now consider two Prym curves $\xi_1$ and $\xi_2$ supported
respectively on $X_1$ and $X_2$. Suppose that $\sigma\colon
X_1\stackrel{\sim}{\to} X_2$ induces an isomorphism between $\xi_1$
and $\xi_2$. Let $\overline{\sigma}\colon Z_1\stackrel{\sim}{\to}Z_2$ be
the induced isomorphism between the stable models. Then it is easy to
see that $\overline{\sigma}$ is an isomorphism between the admissible
covers $\Phi(\xi_1)$ and $\Phi(\xi_2)$.
Moreover, any isomorphism between $\Phi(\xi_1)$ and $\Phi(\xi_2)$ is
obtained in this way.
Hence we have an exact sequence of automorphism groups:
\begin{equation}
\label{auto}
1\rightarrow \Aut_0(\xi)\rightarrow\Aut(\xi)
\rightarrow \Aut(C_{\xi}\rightarrow Z)\rightarrow 1.
\end{equation}

We show that $\Phi$ is surjective.
Let $C\rightarrow Z$ be an admissible double cover,
$I\subset C$ the set of 
fixed points of the involution and $J\subset Z$ their
images. Let $\w{C}\rightarrow C$ and $\nu\colon\w{X}\rightarrow Z$ be the
normalizations of $C$ at $I$ and of $Z$ at $J$ respectively. 
Then $i$ extends to an
involution on $\w{C}$, whose quotient is $\w{X}$, namely:
$\w{C}\rightarrow\w{X}$ is a double cover, ramified over
$q^1_1,q^2_1,\dotsc,q^1_r,q^2_r$, where $r=|J|$ and 
$\nu(q^1_i)=\nu(q^2_i)\in J$
for $i=1,\dotsc,r$.
  Let $L\in\Pic\w{X}$ be the associated line bundle,
satisfying $L^{\otimes 2}\simeq 
\mathcal{O}_{\w{X}}(q^1_1+q^2_1+ \cdots +q^1_r+q^2_r
)$. Finally let $X$ be the quasistable curve
obtained by attaching to $\w{X}$ $r$ rational components
$E_1,\dotsc,E_r$ such that $E_i\cap \w{X}=\{q^1_i,q^2_i\}$.
 Choose  $\eta\in\Pic X$ having degree 1 on all $E_i$ and such that 
$\eta\vert_{\w{X}}= L^{\otimes(-1)}$.
Let $\beta\colon \eta^{\otimes 2}\to
\ol{X}$ be a homomorphism which agrees with  $\eta\vert_{\w{X}}\simeq 
\mathcal{O}_{\w{X}}(-q^1_1-q^2_1- \cdots -q^1_r-q^2_r
)\hookrightarrow\ol{X}$ on $\w{X}$. Then
$\xi=(X,\eta,\beta)$ is a Prym curve with $\eta\not\simeq\mathcal{O}_X$,
and  $C\to Z$ is $\Phi(\xi)$. For
different
choices of $\eta$, the corresponding Prym curves differ by an
inessential isomorphism.
\begin{Proposition}
The map $\Phi$ just defined induces an isomorphism
$$\widehat{\Phi}\colon\prymbar^+\longrightarrow\rg$$
over $\mgbar$.
\end{Proposition}
\begin{proof}
By what precedes, 
$\Phi$ induces a bijection
$\widehat{\Phi}\colon\prymbar^+\to\rg$. 
The statement will follow if  we prove that 
$\widehat{\Phi}$ is a local
isomorphism at every point of $\prymbar^+$. 

Fix a point $\xi=(X,\eta,\beta)\in \prymbar^+$ and consider its universal
deformation
 $(\mathcal{X} \to B, \boldsymbol{\eta},\boldsymbol{\beta})$ 
constructed in 1.2.
Keeping the notations of 1.2,
the line bundle $\boldsymbol{\eta}^{\otimes(-1)}$ determines a double
cover $\overline{\mathcal{P}}\rightarrow\mathcal{X}$, ramified over
$\mathcal{E}_1,\dotsc,\mathcal{E}_r$. 
The divisor $\mathcal{E}_i$ is a $\pr{1}$-bundle over $V_i\subset B$,
and the restriction of its normal bundle to a non trivial fiber
$F$ is
$(\mathcal{N}_{\mathcal{E}_i/\mathcal{X}})_{|F}
\simeq\mathcal{O}_{\pr{1}}(-2)$. 
The inverse image 
$\overline{\mathcal{E}}_i$
of $\mathcal{E}_i$ in $\overline{\mathcal{P}}$ is again a
$\pr{1}$-bundle over $V_i\subset B$, but now 
 the restriction of its normal bundle to a non trivial fiber
$\overline{F}$ is
$(\mathcal{N}_{\overline{\mathcal{E}}_i/
\overline{\mathcal{P}}})_{|\overline{F}}\simeq
\mathcal{O}_{\pr{1}}(-1)$. 
Let $\overline{\mathcal{P}}\rightarrow\mathcal{P}$ be the blow-down 
of $\overline{\mathcal{E}}_1,\dotsc,\overline{\mathcal{E}}_r$. 
We get a diagram
\begin{equation*}
\xymatrix{{\overline{\mathcal{P}}}\ar[d]\ar[r]&{\mathcal{P}}
\ar[d]&
\\{\mathcal{X}}\ar[r]&{\mathcal{Z}}\ar[r]&{B}}
\end{equation*}
where $\mathcal{P}\to \mathcal{Z}\rightarrow B$ is a family
of admissible double covers whose central fiber 
is $C_{\xi}\to Z$. Therefore, up to shrinking $B$, 
there exists a morphism $B\rightarrow W$
such that $\mathcal{P}\to \mathcal{Z}\rightarrow B$ is obtained by
pull-back from the universal deformation
$\mathcal{C}\to\mathcal{Q}\to W$ of $C_{\xi}\to Z$.
Now notice that $\mathcal{Q}\to W$ is a family of stable curves of 
genus $g$, with $Z$ as central fiber: so (again up to shrinking)
it must be a pull-back of 
the universal deformation
$\mathcal{Z}'\to B'$. In the end we get a diagram:
\begin{equation*}
\xymatrix{{\overline{\mathcal{P}}}\ar[d]\ar[r]&{\mathcal{P}}
\ar[d]\ar[r]&{\mathcal{C}}\ar[d]&\\
{\mathcal{X}}\ar[r]\ar[dr]&{\mathcal{Z}}\ar[d]\ar[r]&
{\mathcal{Q}}\ar[d]\ar[r]&{\mathcal{Z}'}\ar[d]\\
&B\ar[r]^{\ph}\ar@(dl,dr)[]_(.4){\Aut (\xi)}
&W\ar[r]^{\psi}\ar@(dl,dr)[]_{\Aut (C_{\xi}\to Z)}
&{B'}\ar@(dl,dr)[]_(.6){\Aut (Z)}
}
\end{equation*}

\vspace{-125pt}\hfill {\small(genus $2g-1$)}

\vspace{25pt}

\hfill {\small(genus $g$)}

\vspace{25pt}

\hfill {\small(bases)}

\vspace{40pt}
\noindent We can assume that $\ph$ and $\psi$ are
surjective. Observe that both maps are equivariant with respect to the
actions of the automorphism groups indicated in the diagram.

Clearly $\ph$ is just the restriction of $\Phi$ to the set of Prym
curves parametrized by $B$.

Now by \eqref{auto},
 $\ph(b_1)=\ph(b_2)$ 
if and only if there exists an inessential
isomorphism between 
$({X}_{b_1},{\eta}_{b_1},\beta_{b_1})$ and
$({X}_{b_2},{\eta}_{b_2},\beta_{b_2})$.
Hence $\ph$ induces an equivariant isomorphism $\widehat{\ph}$:
\begin{equation*}
\xymatrix{
{B/\Aut_0(\xi)}\ar[r]^(.65){\widehat{\ph}}\ar@(dl,dr)[]_(.4){\Aut
  (\xi)/\Aut_0(\xi)} 
&W\ar[r]\ar@(dl,dr)[]_{\Aut (C_{\xi}\to Z)}
&{B'}\ar@(dl,dr)[]_(.6){\Aut (Z)}
}
\end{equation*}
and finally if we mod out by all the automorphism groups, we get
\begin{equation*}
\xymatrix{
{B/\Aut(\xi)}\ar @{_{(}->}[d]
\ar[r]^-{\sim}
&{W/\Aut(C_{\xi}\to Z)}
\ar @{_{(}->}[d]\ar[r]&{B'/\Aut Z} \ar @{_{(}->}[d]\\
{\prymbar^+}\ar[r]^{\widehat{\Phi}}&{\rg}&
{\mgbar.}
}
\end{equation*}
This shows that $\widehat{\Phi}$ is a local isomorphism 
in $\xi$.
\end{proof}
%%%%%%%%%%%%%%%%%%%%%%%%%%%%%%%%%%%%%%%%%%%%%%%%%%%%%%%%%%%%%%%%%%%%%%%%%%
\section{Embedding $\prymbar$ in the compactified Picard variety}\label{geo}
Let $g \ge 3$. For every integer $d$, there is a  
\emph{universal Picard variety}
$$
\pdg \longrightarrow \mgo
$$
whose fiber $J^d(X)$ over a point $X$ of $\mgo$ parametrizes line 
bundles on $X$ of degree $d$, modulo isomorphism.
Denote by $\pry_g^{\,0}$  the inverse image of $\mgo$ under the
finite morphism $\prymbar\to\mgbar$; then we have a commutative diagram
\begin{equation*}
%\label{diagram}
\xymatrix{{\pry_g^{\,0}\ }\ar[dr]\ar@{^{(}->}[r]&{P_{0,g}}\ar[d]\\&{\mgo}}
\end{equation*}
Assume $d \ge 20(g-1)$; this is not a real 
restriction, since for all $t\in\Z_{\geq 0}$ there is a
natural isomorphism $P_{d,g} \cong P_{d+t(2g-2),g}$. Then
$\pdg$ has a natural compactification
$\pdgbar$, endowed with a natural morphism
$\phi_d\colon \pdgbar \to \mgbar$,  
such that $\phi_d^{-1}(\mgo)=\pdg$.
It was constructed in \cite{caporaso} as a GIT quotient
$$
\pi_d\colon H_d \longrightarrow H_d \,/\!/\, G = \pdgbar,
$$
where $G = SL(d-g+1)$ and 
\begin{multline*}
H_d := \{ h \in \Hilb^{dx-g+1}_{d-g}\,|\,h \textrm{ is
 $G$-semistable}\\
\textrm{and the
  corresponding curve is connected }\}
\end{multline*}
(the action of $G$ is linearized by a suitable embedding of
$\Hilb^{dx-g+1}_{d-g}$ in a Grassmannian).

Fix now and in the sequel an integer $t\geq 10$ and define
\begin{align*}
K_{2t(g-1)} := \{& h \in \Hilb^{2t(g-1)x-g+1}_{2t(g-1)-g}| 
        \textrm{ there is a Prym curve } (X, \eta,\beta) 
 \textrm{ and}\\ &\textrm{an embedding }
 h_t\colon X \to \P^{2t(g-1)-g} 
 \textrm{ induced by }
 \eta \otimes \omega_X^{\otimes t},\\
 &\textrm{such that $h$ is the Hilbert point of }  h_t(X) \}.
\end{align*}
Our result is the following.
\begin{Theorem}\label{injection}
The set $K_{2t(g-1)}$ is contained in $H_{2t(g-1)}$; consider its projection 
\[\varPi_t := \pi_{2t(g-1)}(K_{2t(g-1)})\subset\overline{P}_{2t(g-1),g}.\]
There is a natural injective morphism
$$
f_t\colon \prymbar \longrightarrow \overline{P}_{2t(g-1),g}
$$
whose image is $\varPi_t$.
\end{Theorem}
\noindent In particular, the Theorem implies that $\varPi_t$ is a closed
subvariety of $\overline{P}_{2t(g-1),g}$.

The proof of Theorem \ref{injection} 
will be achieved in several steps and will take  the rest of this
section. The argument is the one used in \cite{fontanari} to show
the existence of an injective morphism
$\overline{S_g}\rightarrow\overline{P}_{(2t+1)(g-1),g}$ of the
moduli space of spin curves in the corresponding compactified Picard
variety.

One can define (see \cite{caporaso}, \S 8.1) the 
contravariant functor $\pdgfun$ from the category of
schemes to the one of sets, which to every scheme 
$S$ associates the set $\pdgfun(S)$ of equivalence classes of polarized 
families of quasistable curves of genus $g$ 
$$
f\colon (\mathcal{X}, \mathcal{L}) \longrightarrow S
$$
such that $\mathcal{L}$ is a relatively very ample line bundle of degree 
$d$ whose multidegree satisfies the following Basic Inequality on each 
fiber. 
\begin{Definition}
Let $X = \bigcup_{i=1}^n X_i$ be a projective, nodal, connected curve of 
arithmetic genus $g$, where the $X_i$'s are the irreducible components of $X$. 
We say that the multidegree $(d_1, \dotsc, d_n)$ 
satisfies the \emph{Basic Inequality} 
if for every complete subcurve $Y$ of $X$ 
of arithmetic genus $g_Y$ we have
$$
m_Y \le d_Y \le m_Y + k_Y
$$
where
$$
d_Y = \sum_{X_i \subseteq Y} d_i,\quad
k_Y= \vert Y \cap \overline{X \smallsetminus Y} \vert \quad\text{and}\quad
m_Y= \frac{d}{g-1} \left( g_Y-1+ \frac{k_Y}{2} \right) - \frac{k_Y}{2}
$$
(see \cite{caporaso} p.~611 and p.~614).
\end{Definition}
\noindent Two families over $S$, $(\mathcal{X}, \mathcal{L})$ and 
$(\mathcal{X}', \mathcal{L}')$ are equivalent if there exists 
an $S$-isomorphism $\sigma\colon \mathcal{X} \to \mathcal{X}'$
and a line bundle $M$ on $S$ such that
$\sigma^*\mathcal{L}' \cong \mathcal{L} \otimes f^*M$.

By \cite{caporaso}, Proposition~8.1, there is a morphism of functors:
\begin{equation}\label{morphism}
\pdgfun \longrightarrow \Hom (\hspace{0.1cm} \cdot \hspace{0.1cm}, \pdgbar)
\end{equation}
and $\pdgbar$ coarsely represents $\pdgfun$ if and only if 
\begin{equation}\label{coprime}
(d-g+1, 2g-2)=1.
\end{equation}
\begin{Proposition}\label{naturalmorphism}
For every integer $t \ge 10$ there is a natural morphism:
$$
f_t\colon \prymbar \longrightarrow \overline{P}_{2t(g-1),g}.
$$
\end{Proposition}
\begin{proof} 
First of all, notice that in this case (\ref{coprime}) does not hold, 
so the points of $\overline{P}_{2t(g-1),g}$ 
are \emph{not} in one-to-one correspondence with 
equivalence classes of very ample line bundles of degree $2t(g-1)$ 
on quasistable 
curves, satisfying the Basic Inequality (see \cite{caporaso}, p.~654). 
However, we claim that the thesis can be deduced from the existence of a 
morphism of functors:
\begin{equation}\label{functor}
F_t\colon \prymfun \longrightarrow \overline{\mathcal{P}}_{2t(g-1),g}.
\end{equation}
Indeed, since $\prymbar$ coarsely represents $\prymfun$, any morphism of 
functors $\prymfun \to \Hom (\hspace{0.1cm} \cdot \hspace{0.1cm}, T)$
induces a morphism of schemes $\prymbar \to T$, so the claim follows from 
(\ref{morphism}).
Now, a morphism of functors as (\ref{functor}) is the datum for any scheme 
$S$ of a set-theoretical map 
$$
F_t(S)\colon \prymfun(S) \longrightarrow \overline{\mathcal{P}}_{2t(g-1),g}(S),
$$
satisfying obvious compatibility conditions.
Let us define $F_t(S)$ in the following way:
$$
( f\colon \mathcal{X} \to S, \boldsymbol{\eta},\boldsymbol{\beta}) 
\ \mapsto\ ( f\colon (\mathcal{X}, \boldsymbol{\eta} 
\otimes \omega_f^{\otimes t}) 
\to S ).
$$
In order to prove that $F_t(S)$ is well-defined, the only non-trivial
matter is to check that the multidegree of $\boldsymbol{\eta}
\otimes \omega_f^{\otimes 
t}$ satisfies the Basic Inequality on each fiber, so the thesis follows 
from the next Lemma.
\end{proof}
\begin{Lemma}\label{basic}
Let $(X,\eta,\beta)$ be a Prym curve.
If $Y$ is a complete subcurve of $X$ and $d_Y$ is the degree of 
$(\eta \otimes \omega_X^{\otimes t})\vert_Y$, then 
$m_Y \le d_Y \le m_Y + k_Y$ in the notation of the Basic Inequality. 
Moreover, if $d_Y = m_Y$ then
$\w{k}_Y := \vert \w{Y} \cap \overline{\w{X}
  \smallsetminus \w{Y}}  
\vert = 0$. 
\end{Lemma}
\begin{proof} In the present case, the Basic Inequality simplifies as follows:
$$
- \frac{k_Y}{2} \le e_Y \le \frac{k_Y}{2},
$$ 
where $e_Y := \deg \eta \vert_Y$.
By the definition of a Prym curve, the degree $e_Y$ depends only on 
the exceptional components of $X$ intersecting $Y$. 

\noindent For any exceptional component $E$ of $X$ with $E \subseteq 
\overline{X \smallsetminus Y}$, let $m := \vert E \cap Y \vert$.  
The contribution of $E$ to $k_Y$ is $m$, while its contribution to 
$e_Y$ is $- \frac{m}{2}$.

\noindent Next, for any exceptional component $E$ of $X$ with $E \subseteq Y$, 
let $l := \vert E \cap \overline{X \smallsetminus Y} \vert$.  
The contribution of $E$ to $k_Y$ is $l$, while its contribution to 
$e_Y$ is $1 - \frac{2 - l}{2}=\frac{l}{2}$.

\noindent Summing up, we see that the Basic Inequality holds. Finally, 
if $\w{k}_Y \ne 0$, then there exists a non-exceptional component 
of $X$ intersecting $Y$. Such a component contributes at least $1$ to 
$k_Y$, but it does not affect $e_Y$; hence $- \frac{k_Y}{2} < e_Y$ 
and the proof is over.
\end{proof}
By applying \cite{caporaso}, Proposition~6.1, from the first part of 
Lemma~\ref{basic} we deduce 
$$
K_{2t(g-1)} \subset H_{2t(g-1)}.
$$
Moreover, the second part of the same Lemma provides a crucial 
information on Hilbert points corresponding to Prym curves.
\begin{Lemma}\label{closed}
If $h\in K_{2t(g-1)}$,
then the orbit of $h$ is closed in the semistable locus. 
\end{Lemma}
\begin{proof} 
Let $(X,\eta,\beta)$ be a Prym curve such that $h$ is the Hilbert point of
an embedding $h_t\colon X\to \P^{2t(g-1)-g}$ induced by 
 $\eta \otimes \omega_X^{\otimes t}$.
Just recall the first part of \cite{caporaso}, Lemma~6.1, which 
says that the orbit of $h$ is closed in the semistable locus if
and only  
if $\w{k}_Y=0$ for every subcurve $Y$ of $X$ such that $d_Y = m_Y$, so 
the thesis is a direct consequence of Lemma~\ref{basic}.  
\end{proof}
\begin{proof}[Proof of Theorem \ref{injection}]
It is easy to check that $f_t(\prymbar)=\varPi_t$. Indeed,
if $(X, \eta,\beta) \in \prymbar$, then any choice of a base 
for $H^0(X, \eta \otimes \omega_X^{\otimes t})$ gives an 
embedding $h_t\colon X \to \P^{2t(g-1)-g}$  and 
$f_t(X, \eta,\beta)
= \pi_{2t(g-1)}(h)$, where $h\in K_{2t(g-1)}$ is the
Hilbert point of $h_t(X)$.
Conversely, if $\pi_{2t(g-1)}(h) \in \varPi_t$, then there is a 
Prym curve $(X, \eta,\beta)$ and an embedding 
$h_t\colon X \to \P^{2t(g-1)-g}$ such that 
$h$ is the Hilbert point of 
$h_t(X)$ and $f_t(X, \eta,\beta)= \pi_{2t(g-1)}(h)$.

Next we claim that $f_t$ is injective. Indeed, let $(X, \eta,\beta)$ 
and $(X', \eta',\beta')$ be two Prym curves and assume that 
$f_t(X, \eta,\beta) = f_t(X', \eta',\beta')$.
Choose bases for $H^0(X, \eta \otimes \omega_X^{\otimes t})$ and
$H^0(X', \eta' \otimes \omega_{X'}^{\otimes t})$ and 
embed $X$ and $X'$ in $\P^{2t(g-1)-g}$.
If $h$ and $h'$ are the corresponding Hilbert points, 
then $\pi_{2t(g-1)}(h)= \pi_{2t(g-1)}(h')$ and
the Fundamental Theorem of GIT implies that 
$\overline{O_G(h)}$ and $\overline{O_G(h')}$ intersect in the 
semistable locus. It follows from Lemma~\ref{closed} that 
$O_G(h) \cap O_G(h') \ne \emptyset$, so
$O_G(h)=O_G(h')$ and there is an isomorphism
$\sigma\colon (X,\eta,\beta) \to (X',\eta',\beta')$. 
\end{proof}
Observe that Theorem \ref{injection} and Lemma \ref{closed} 
imply that $K_{2t(g-1)}$
is a constructible set in $H_{2t(g-1)}$.

%%%%%%%%%%%%%%%%%%%%%%%%%%%%%%%%%%%%%%%%%%%%%%%%%%%%%%%%%%%%%%%%%%%
\section{Fiberwise description}\label{com}
Let $Z$ be a stable curve of genus $g$. We recall that the dual graph
$\Gamma_Z$ of $Z$ is the graph whose vertices are the
irreducible components of $Z$ and whose edges are the nodes of $Z$. 
The first Betti number of $\Gamma_Z$ is
$b_1(\Gamma_Z)=\delta-\gamma+1=g-g^{\nu}$, where $\delta$ is the number
of nodes of $Z$, $\gamma$ the number of its irreducible components and
$g^{\nu}$ the genus of its normalization.

We denote by $\prymx$ the scheme parametrizing Prym curves $(X,\eta,\beta)$ 
such that the stable model of $X$ is $Z$, 
modulo inessential isomorphisms, 
and by $S_Z$ the 
analogue for spin curves. Since by Lemma \ref{inessential} the homomorphism
$\beta$ is not relevant in determining the inessential isomorphism
class of $(X,\eta,\beta)$, 
in this section we will omit it and just write $(X,\eta)$.

When $\Aut(Z)=\{\text{Id}_Z \}$, $\prymx$ is the scheme-theoretical
fiber over
$[Z]$ of the morphism $p\colon\prymbar\rightarrow \mgbar$.
Recall that $p$
is finite of degree $2^{2g}$, and \'etale over $\mgo$. 

For any 0-dimensional scheme $P$ we denote by $L(P)$ the set of
integers occurring as
multiplicities of components of $P$.

In this section we describe the numerical properties of $\prymx$, namely
the number of irreducible components and their multiplicities, 
showing that they depend
only on the dual graph $\Gamma_Z$ of
$Z$. Using this, we give some properties of $L(\prymx)$, and show that
in some cases the set of multiplicities
$L(\prymx)$ gives informations on $Z$. In particular, we
show that the morphism $\prymbar\rightarrow \mgbar$ is \'etale 
over $\mgbar^0\smallsetminus D_{irr}$.

We use the techniques and results of~\cite{spin}, where the same
questions about the numerics of $S_Z$ are studied (see also \cite{cs}, \S~3). 
Quite surprisingly, the schemes $P_Z$ and $S_Z$ are not isomorphic in general.

Finally we will show with an example that, differently from the case of 
spin curves, the set of multiplicities $L(\prymx)$
appearing in $\prymx$ does not always identify 
curves having two smooth components.

Let $X$ be a quasistable curve having $Z$ as stable model 
and consider the set
\[ \Delta_X:=\{z\in\Sing Z\,|\, z\text{ is \emph{not} the image of an
  exceptional component of }X\}.\]
%Clearly the complementary set in $\Sing Z$ is
%\[ \Delta_X^c=\{z\in\Sing Z\,|\, z\text{ is the image of an
%  exceptional component of }X\}.\]
Given $Z$, the quasistable curve $X$ is determined by 
$\Delta_X$, or equivalently by 
$$\Delta_X^c:=\Sing Z\smallsetminus\Delta_X
=\{\text{images in $Z$ of the 
 exceptional components of }X\}. $$ 
Remark that
any subset of $\Sing Z$  can be seen as
a subgraph of the dual graph $\Gamma_Z$ of $Z$.
%Let $\w{X}=\overline{X\smallsetminus(\cup_iE_i)}$, where $E_i$ are the
%exceptional components of $X$. Remark that the dual graph of $\w{X}$ is
%exactly $\Delta_X$. 

We recall that the \emph{valency} of a vertex of a graph
is the number of edges 
ending in that vertex and a graph $\Gamma$ is \emph{eulerian} if it has 
all even valencies. Thus $\Gamma_Z$ is eulerian if and only if for any
irreducible component $C$ of $Z$, $|C\cap \overline{Z\smallsetminus
  C}|$ is even. 
The set  $\mathcal{C}_{\Gamma}$ of all eulerian
subgraphs of $\Gamma$ is called the cycle space of  $\Gamma$. 
There is a natural identification of
$\mathcal{C}_{\Gamma}$ with $H_1(\Gamma,\Z_2)$, so
$|\mathcal{C}_{\Gamma}|= 2^{b_1(\Gamma)}$ (see \cite{spin}). 
Reasoning as in \cite{spin}, Section 1.3, we can show the following:
\begin{Proposition}
\label{conti}
Let $X$ be a quasistable curve having $Z$ as stable model.

The curve $X$ is the support of a Prym curve  if and only if
$\Delta_X^c$ is eulerian.

If so, there are
$2^{2g^{\nu}+b_1(\Delta_X)}$ different choices for $\eta\in\Pic X$ such that
$(X,\eta)\in\prymx$. 

For each such $\eta$, the point $(X,\eta)$ has multiplicity
$2^{b_1(\Gamma_Z)-b_1(\Delta_X)}$ in $\prymx$.
\end{Proposition}
\noindent
Hence the number of irreducible components of $\prymx$ is
\[
2^{2g^{\nu}}\cdot\sum_{\Sigma\in\mathcal{C}_{\Gamma_Z}}2^{b_1(\Sigma^c)},\]
and its set of multiplicities is given by
\[ L(\prymx)=\{ 2^{b_1(\Gamma_Z)-b_1(\Delta)}\,|\,\Delta^c\in
\mathcal{C}_{\Gamma_Z}\}. \]
 Remark 
that since $|\mathcal{C}_{\Gamma_Z}|=2^{b_1(\Gamma_Z)}$, we can
check immediately from the proposition that the length of $\prymx$ is 
\[ \sum_{\Sigma\in\mathcal{C}_{\Gamma_Z}}(2^{2g^{\nu}+b_1(\Sigma^c)}
\cdot 2^{b_1(\Gamma_Z)-b_1(\Sigma^c)})=2^{b_1(\Gamma_Z)}\cdot 
2^{2g^{\nu}+b_1(\Gamma_Z)}=2^{2g}.\]
As a consequence of Proposition \ref{conti}, 
we see that
\begin{enumerate}[$\bullet$]
\item
{\em a point $(X,\eta)$ in $\prymx$ is non reduced if and only if $X$ is
non stable. }
\end{enumerate}
\begin{Example}[curves having two smooth components]
Let $Z=C_1\cup C_2$, $C_i$ smooth irreducible, 
$|C_1\cap C_2|=\delta\geq 2$.
\begin{center}
  
\scalebox{0.35}{\includegraphics{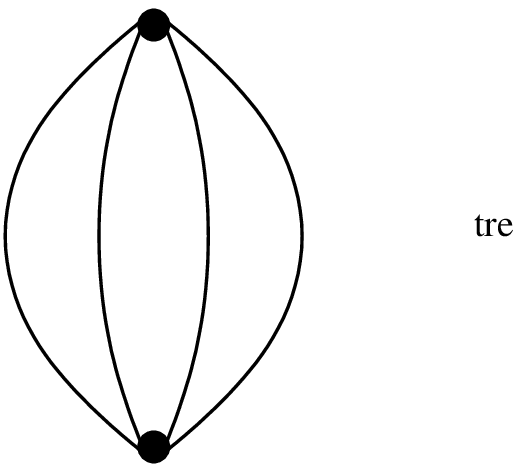}}

\end{center}
Let $X$ be a quasistable curve having $Z$ 
as stable model and let $\Delta_X$ 
be the corresponding subset of $\Sing Z$.
The subgraph $\Delta_X^c$  is eulerian if and 
only if $|\Delta_X^c|$ is even. Therefore $X$
 is support of a Prym curve if and only if it has an even
number $2r$ of exceptional components. If so, for each choice of
$\eta\in\Pic X$ such that $(X,\eta)\in\prymx$, this point will have
multiplicity $2^{b_1(\Gamma_Z)-b_1(\Delta_X)}$. We have
$b_1(\Gamma_Z)=\delta-1$ and $|\Delta_X|=\delta-2r$, so
\[
b_1(\Delta_X)=\begin{cases} \delta-2r-1\quad &\text{if }
\ 2r\leq\delta-2,\cr
0\quad &\text{if }\ \delta-1\leq 2r\leq \delta
\end{cases}
\]
and we get
\[ L(\prymx)=\{2^{2r}\,|\,0\leq r\leq \frac{1}{2}
\delta -1\}\cup\{2^{\delta-1}\}.\]
\end{Example}
\begin{Proposition}[combinatorial properties of $L(\prymx)$]
\label{combprop}
The following properties hold:
\begin{enumerate}[(1)]
\item $1\in L(\prymx)$;
\item $\max L(\prymx)=2^{b_1(\Gamma_Z)}$;
\item $2^g\in L(\prymx)$ if and only if $Z$ has only rational components;
\item $\prymx$ is reduced if and only if $Z$ is of compact type;
\item if $\Gamma_Z$ is an eulerian
  graph, then
$L(\prymx)=L(S_Z)$.
\end{enumerate}
\end{Proposition}
\begin{proof}
(1) Choosing $\Delta_X=\Gamma_Z$, we get $X=Z$;
since the empty set is trivially in $\mathcal{C}_{\Gamma_Z}$,
there always exists $\eta\in \Pic Z$ such that
$(Z,\eta)\in\prymx$. This $\eta$ is a square root of $\ol{Z}$; there are 
$2^{2g^{\nu}+b_1(\Gamma_Z)}$ choices for it, and it will appear with
multiplicity 1 in $\prymx$. So $1\in L(\prymx)$.

(2) From Proposition \ref{conti} we get $\max  L(\prymx)\leq
    2^{b_1(\Gamma_Z)}$. 
Set $M=\max\{b_1(\Sigma)\,|\,\Sigma\in\mathcal{C}_{\Gamma_Z}\}$ and
let $\Sigma_0 \in\mathcal{C}_{\Gamma_Z}$ be such that
$b_1(\Sigma_0)=M$. 
By Proposition \ref{conti}, we know that
$2^{b_1(\Gamma_Z)-b_1(\Sigma_0^c)}
\in L(\prymx)$.
We claim that $b_1(\Sigma_0^c)=0$. Indeed, if not,
$\Sigma_0^c$ contains a subgraph $\sigma$ with $b_1(\sigma)=1$
and having all valencies equal to 2. Then 
$\Sigma_0\cup\sigma\in\mathcal{C}_{\Gamma_Z}$ and
$b_1(\Sigma_0\cup\sigma)>M$, a contradiction.
Hence we have points of multiplicity
$2^{b_1(\Gamma_Z)}$ in $\prymx$, so $\max  L(\prymx)=2^{b_1(\Gamma_Z)}$.

Property (3) is immediate from (2), since $b_1(\Gamma_Z)=g$ if and
only $g^{\nu}=0$.

Also property (4) is immediate from (2), because $L(\prymx)=\{1\}$ if
and only if $b_1(\Gamma_Z)=0$.

(5) Assume that $\Gamma_Z$ is eulerian. Then
$\Delta_X^c\in\mathcal{C}_{\Gamma_Z}$ if and only if
$\Delta_X\in\mathcal{C}_{\Gamma_Z}$, so we have
\[ L(\prymx)=L(S_Z)=\{ 2^{b_1(\Gamma_Z)-b_1(\Delta_X)}\,|\,\Delta_X\in
\mathcal{C}_{\Gamma_Z}\}\]
(see \cite{spin} for the description of $L(S_Z)$).
\end{proof}
Property (4) implies the following
\begin{Corollary}
The morphism $p\colon\prymbar\to\mgbar$ is \'etale over
$\mgbar^0\smallsetminus D_{irr}$.
\end{Corollary} 
Consider now property (1) of Proposition \ref{combprop}.
It shows, in particular, that in general $\prymx$ and
$S_Z$ are not 
isomorphic and do not have the same set of multiplicities: indeed, for
spin curves, it can very well happen that $1\not\in L(S_Z)$ (see
example after Corollary \ref{cor}).

The following shows that in some cases, $L(\prymx)$ gives informations
on $Z$.
\begin{Corollary}
\label{cor}
Let $Z$ be a stable curve and $\nu\colon Z^{\nu}\rightarrow Z$ its
normalization. Assume that for every irreducible component
$C$ of $Z$, the number $|\nu^{-1}(C\cap\Sing Z)|$ 
 is even
and at least $4$. 
\begin{enumerate}[(i)]
\item If $2^{b_1(\Gamma_Z)-2}\not\in
  L(\prymx)$, then $Z=C_1\cup C_2$, with
 $C_1$ and $C_2$ smooth and irreducible.
\begin{center}
\scalebox{0.30}{\includegraphics{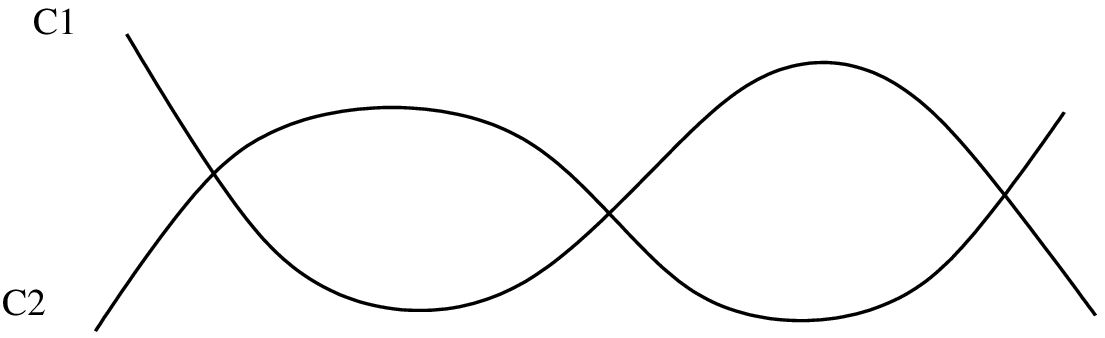}}
\end{center}
\item If $2^{b_1(\Gamma_Z)-3}\not\in
  L(\prymx)$, then either $Z$ is irreducible with two nodes, or 
$Z=C_1\cup C_2\cup C_3$, with $C_i$ smooth
  irreducible and $|C_i\cap C_j|=2$ for $1\leq i<j\leq 3$.
\begin{center}
\scalebox{0.30}{\includegraphics{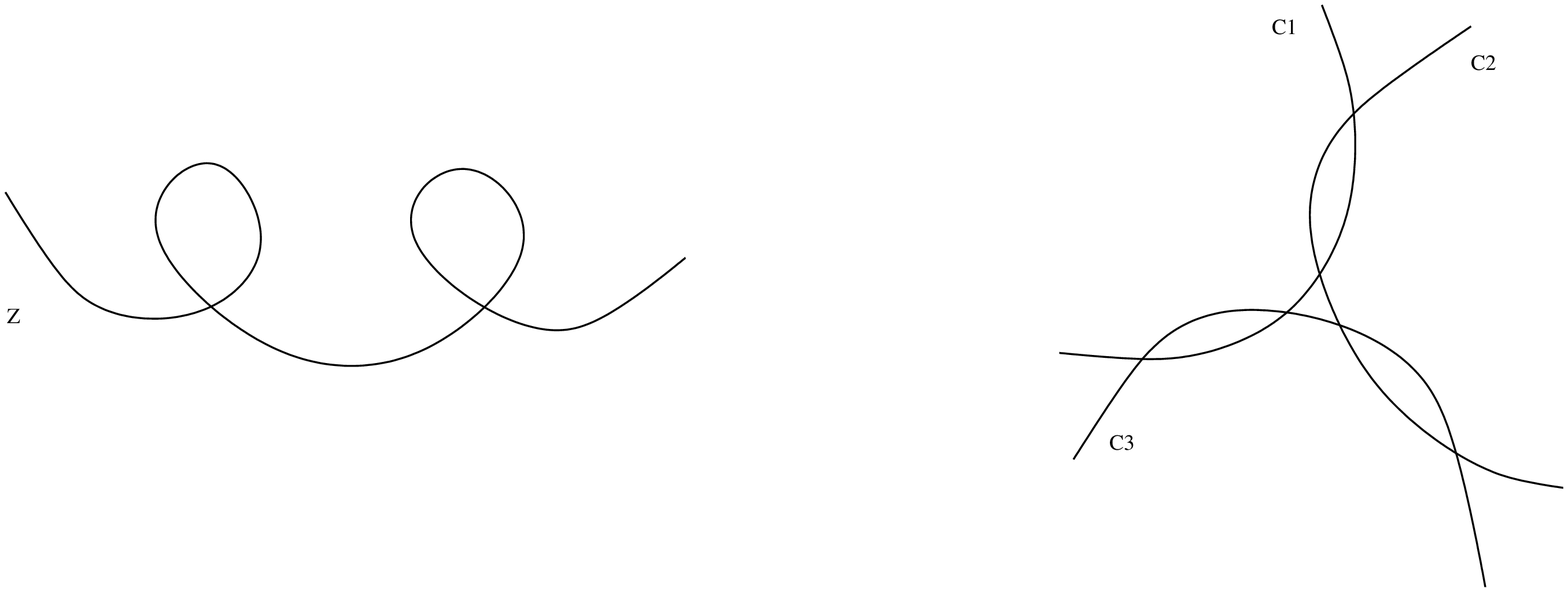}}
\end{center}
\end{enumerate}
\end{Corollary}
\begin{proof}
By hypothesis $\Gamma_Z$ is eulerian, so property (5) says that
$L(\prymx)=L(S_Z)$. Then $(i)$ follow immediately from \cite{spin}, 
Theorem~11. Let us show $(ii)$. If $b_1(\Gamma_Z)\geq 4$, 
by property (2) we can apply \cite{spin}, Theorem~13; then 
$Z$ has three smooth components meeting each
other in two points. Assume $b_1(\Gamma_Z)\leq 3$ and let $\delta$,
$\gamma$ be the number of nodes and of irreducible components of $Z$.
Since all vertices of $\Gamma_Z$ have valency at least 4, we have
$\delta\geq 2\gamma$, so $\gamma\leq b_1(\Gamma_Z)-1\leq 2$. Then by
an easy check we see that the only possibility which satisfies all the
hypotheses is $\gamma=1$ and $\delta=2$. 
\end{proof}
In \cite{spin} it is shown (Theorem 11) that 
$L(S_Z)$  allows to recover curves having two smooth components. 
Instead, when the number of nodes is odd, it is no more true that 
these curves are characterized by $L(\prymx)$. 
For instance, consider the graphs:
\begin{center}
  
\scalebox{0.35}{\includegraphics{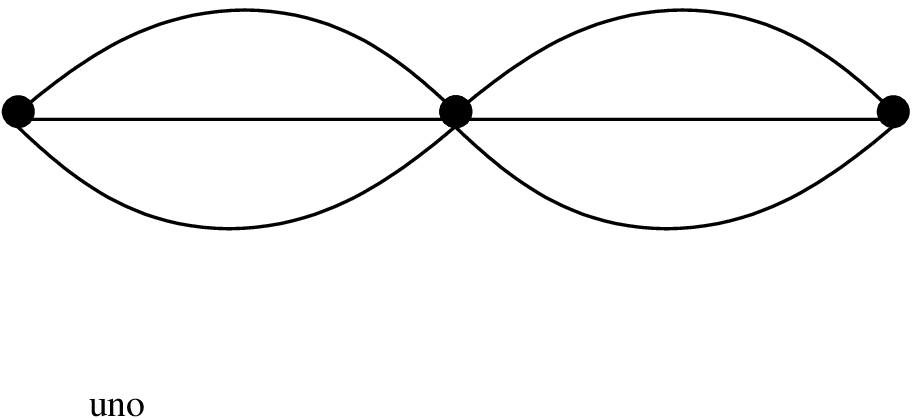}}\hspace{100pt}
\scalebox{0.35}{\includegraphics{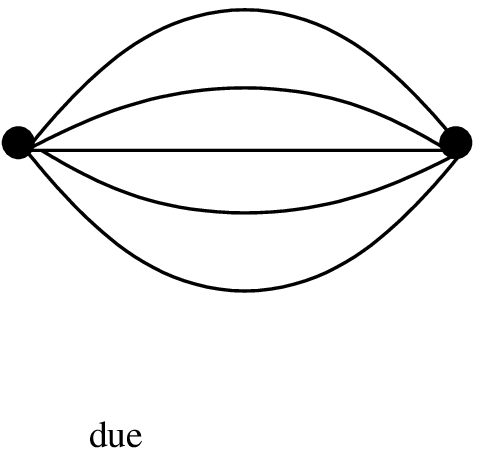}}

\end{center}
It is easy to see that if $Z_1$, $Z_2$ are stable curves with
$\Gamma_{Z_i}=\Gamma_i$ for $i=1,2$, we have
$L(P_{Z_1})=L(P_{Z_2})=\{1,4,16\}$, while $L(S_{Z_1})=\{4,8,16\}$ and
$L(S_{Z_2})=\{2,8,16\}$.
%Even if $\Gamma_Z=\Delta_X\cup\Delta_X^c$, it is false in
%general that $b_1(\Gamma_Z)-b_1(\Delta_X)=b_1(\Delta_X^c)$: we just
%have $b_1(\Gamma_Z)-b_1(\Delta_X)\geq b_1(\Delta_X^c)$.
%Hence, for Prym curves we cannot express the set $L(\prymx)$ 
%in terms of the set of Betti numbers of eulerian subgraphs 
%of $\Gamma_Z$, as it has been done in \cite{spin} for spin curves.
%%%%%%%%%%%%%%%%%%%%%%%%%%%%%%%%%%%%%%%%%%%%%%%%%%%%%%%%%%%%%%%%%%%

\bigskip

\small
%\bibliographystyle{alpha}
%\bibliography{Biblio}

\begin{thebibliography}{ACGH85}

\bibitem[ABH01]{prymhulek}
Valery Alexeev, Christina Birkenhake, and Klaus Hulek.
\newblock Degenerations of {P}rym varieties.
\newblock {\em J.\ Reine Angew.\ Math.},
  553:73--116, 2001.

\bibitem[ACGH85]{acgh}
Enrico Arbarello, Maurizio Cornalba, Phillip Griffiths, and Joseph Harris.
\newblock {\em Geometry of Algebraic Curves, I}, volume 267 of {\em Grundlehren
  der mathematischen Wissenschaften}.
\newblock Springer-Verlag, 1985.

\bibitem[Bea77]{beauville}
Arnaud Beauville.
\newblock {P}rym varieties and the {S}chottky problem.
\newblock {\em Invent.\ Math.}, 41:149--196, 1977.

\bibitem[Bea96]{beau}
Arnaud Beauville.
\newblock {\em Complex Algebraic Surfaces}.
\newblock Cambridge University Press, second edition, 1996.

\bibitem[Ber99]{mira}
Mira Bernstein.
\newblock {\em Moduli of curves with level structure}.
\newblock PhD thesis, Harvard University, 1999.
\newblock Available at the web page \verb+http://www.math.berkeley.edu/~mira+.

\bibitem[Cap94]{caporaso}
Lucia Caporaso.
\newblock A compactification of the universal {P}icard variety over the moduli
  space of stable curves.
\newblock {\em J.\ Amer.\ Math.\ Soc.}, 7(3):589--660,
  1994.

\bibitem[CC03]{spin}
Lucia Caporaso and Cinzia Casagrande.
\newblock Combinatorial properties of stable spin curves.
\newblock {\em Comm.\ Algebra}, 31(8):3653--3672, 2003.
\newblock Special Issue in Honor of Steven L.\ Kleiman.

\bibitem[Cor89]{cornalba1}
Maurizio Cornalba.
\newblock Moduli of curves and theta-characteristics.
\newblock In {\em Lectures on {R}iemann Surfaces: Proceedings of the {C}ollege
  on {R}iemann Surfaces, {I}nternational {C}entre for {T}heoretical {P}hysics,
  Trieste, 1987}, pages 560--589. World Scientific, 1989.

\bibitem[Cor91]{cornalba2}
Maurizio Cornalba.
\newblock A remark on the {P}icard group of spin moduli space.
\newblock {\em Atti dell'Accademia Nazionale dei Lincei, Matematica e
  Applicazioni}, IX(2):211--217, 1991.

\bibitem[CS03]{cs}
Lucia Caporaso and Edoardo Sernesi.
\newblock Characterizing curves by their odd theta-characteristics.
\newblock {\em J.\ Reine Angew.\ Math.},
  562:101--135, 2003.

\bibitem[DM69]{dm}
Pierre Deligne and David Mumford.
\newblock The irreducibility of the space of curves of given genus.
\newblock {\em Publ.\ Math.\ I.\ H.\
  E.\ S.}, 36:75--109, 1969.

\bibitem[DS81]{prym}
Ron Donagi and Roy~Campbell Smith.
\newblock The structure of the {P}rym map.
\newblock {\em Acta Math.}, 146:25--102, 1981.

\bibitem[Fon02]{fontanari}
Claudio Fontanari.
\newblock On the geometry of the compactification of the universal {P}icard
  variety.
\newblock {Preprint math.AG/0202168}, 2002.

\bibitem[Ray71]{SGA1EXPXII}
Michel Raynaud.
\newblock Expos{\'e} {XII}. {G}{\'e}om{\'e}trie alg{\'e}brique et
  g{\'e}om{\'e}trie analytique.
\newblock In {\em Rev{\^e}tements {\'E}tales et Groupe Fondamental.
  S{\'e}minaire de G{\'e}om{\'e}trie Alg{\'e}brique du Bois Marie 1960--1961
  (SGA1). Dirig{\'e} par Alexandre Grothendieck}, volume 224 of {\em Lecture
  Notes in Mathematics}, pages 311--343. Springer-Verlag, 1971.

\end{thebibliography}

\vspace{0.5cm}

\noindent
Edoardo Ballico and Claudio Fontanari\newline
Universit\`a degli Studi di Trento \newline
Dipartimento di Matematica \newline
Via Sommarive 14 \newline
38050 Povo (Trento) \newline
Italy \newline
e-mail: ballico@science.unitn.it,\ \ 
fontanar@science.unitn.it

\vspace{0.5cm}

\noindent
Cinzia Casagrande \newline
Universit\`a degli Studi Roma Tre \newline
Dipartimento di Matematica \newline
Largo S. L. Murialdo 1 \newline
00146 Roma \newline
Italy \newline
e-mail: casagran@mat.uniroma3.it
\end{document}